\documentclass{amsart}
\usepackage{graphicx}
\theoremstyle{definition}

\theoremstyle{remark}

\numberwithin{equation}{section}

\begin{document}

\title{The advanced maximum principle for parabolic systems on manifolds with boundary}%
\author{Hong Huang}%
\address{School of Mathematics Science,Beijing Normal University,Beijing 100875, P. R. China}%
\email{hhuang@bnu.edu.cn}%

\thanks{Partially supported by NSFC no.10671018.}%


\begin{abstract}
In this short note we extend  Chow and Lu's advanced maximum
principles for parabolic systems on closed manifolds to the case of
compact manifolds with boundary, which also generalizes a Hopf type
theorem of Pulemotov.
\end{abstract} \maketitle

\section {Introduction}

In his important 1982 paper [H1], Hamilton introduced Ricci flow and
used it to prove that any closed three-manifolds with positive Ricci
curvature is diffeomorphic to a spherical space form. In 1986
Hamilton [H2] introduced an advanced maximum principle for parabolic
systems on closed manifolds with evolving metrics, and using it he
was able to prove that any closed 4-manifold with positive curvature
operator is diffeomorphic to a spherical space form, in addition to
giving a simplified proof of his 1982 3-dimensional result.

Later, Hamilton [H3] and Ivey [I] proved independently an important
pinching estimate for 3-dimensional Ricci flow using Hamilton's
advanced maximum principle. For recent applications of Hamilton's
advanced maximum principle see for example B$\ddot{o}$hm and Wilking
[BW], Brendle and Schoen [BS].

In [CL] Chow and Lu presented two useful generalizations (see
Theorems 3 and 4 in [CL]) of Hamilton's advanced maximum principle
for parabolic systems on closed manifolds. In this short note we
extend Chow and Lu's results to the case of compact manifolds with
boundary, which also generalize a Hopf type theorem of Pulemotov
(see Theorem 2.1 in [P],which itself is a generalization of
Hamilton's maximum principle in [H2] and Shen's Hopf type theorem in
[S]).

Below we will follow closely  the notations in [CL]. Let $M$ be a
compact manifold with boundary with a one-parameter family of
Riemannian metrics $g_{ij}(t), 0\le t \le T$ with $T < \infty$.

Let $\pi:V\rightarrow M$ be a vector bundle with a fixed (i.e.
time-independent) bundle metric $h_{ab}$. We equip $V$ with a family
of time-dependent connections $\nabla_t$ compatible with $h_{ab}$.
We define the Laplacian $\Delta_t$ acting on a section $\sigma \in
\Gamma(V)$ by $\Delta_t\sigma=g^{ij}(x,t)(\nabla_t)_i(\nabla_t)_j
\sigma$ as usual.

Let $U$ be an open subset of $V$, and for each $t\in [0,T]$, let
$\mathcal{K}(t)\subset U$ be a closed subset such that in each fiber
$V_x$ over $x\in M$ the subset $\mathcal{K}_x(t)=\mathcal{K}(t)\cap
V_x$ is nonempty, closed and convex, and  $\mathcal{K}(t)$ is
invariant under parallel translation defined by the connection
$\nabla_t$. We  assume further  the space-time track
$\mathcal{T}=\{(v,t)\in V\times [0,T]|v\in \mathcal{K}(t), t\in
[0,T]\}$ is closed.

Let $F:U\times [0,T]\rightarrow V$ be a fiber preserving map,which
may be viewed  as a time-dependent vector field  on  $U$ which is
tangent to the fibers. We assume $F(x,\sigma, t)$ is continuous in
$x, t$ and Lipschitz continuous in $\sigma$ in the sense that the
inequalities

$|F(x,\sigma_1,t)-F(x,\sigma_2,t)|\le C_B|\sigma_1-\sigma_2|$

hold true for all $x\in M,t\in [0,T]$ and $|\sigma_1|\le
B,|\sigma_2|\le B$, where $C_B$ is a constant depending only on $B$.

We will consider sections of $V$ which satisfy the following

(PDE) $\frac{\partial}{\partial t}\sigma(x,t)=\Delta_t
\sigma(x,t)+u^i(x,t)(\nabla_t)_i \sigma(x,t)+F(x,\sigma(x,t),t)$,

where $u^i$ is a time-dependent vector field on $M$.

For convenience we recall some notions about convex sets. Let
$\mathcal{J}$ be a closed convex subset of $R^n$ and $v \in
\partial \mathcal{J}$. We denote by $C_v \mathcal{J}$ the tangent cone to
 $\mathcal{J}$ at $v$, which is the
smallest closed convex cone containing $\mathcal{J}$ with vertex at
$v$, and by $S_v \mathcal{J}$ the set of support functions $l$ for
$\mathcal{J}$ at $v$, which are linear functions on $R^n$ satisfying
$|l|=1 $ and $l(v)\ge l(w)$ for all $w \in \mathcal{J}$.

 Now we can state our

 \hspace *{0.4cm}

{\bf Theorem }  Let $M$,$V$,$\mathcal{K}(t)$ and $F$ be as above.
Assume that for any $x\in M$ and any $t_0\in [0,T)$, any solution
$\rho_x(t)$ of the

(ODE) $\frac{d\rho_x}{dt}=F(x,\rho_x,t)$

which starts in $\mathcal{K}_x(t_0)$ at $t_0$ will remain in
$\mathcal{K}_x(t)$ for all later times $t\in [t_0,T]$. Moreover,
assume that for any $x\in
\partial M$, $t\in (t_0,T)$ (given $t_0\in [0,T)$) and any $v \in
\partial \mathcal{K}_x(t)$, the solution $\sigma$ of the (PDE) satisfies
$v+(\nabla_t)_{\nu}\sigma(x,t)\in C_v\mathcal{K}_x(t)$, where $\nu$
is the unit outward normal vector to $\partial M$ at $(x,t)$. Then
the solution $\sigma$ of the (PDE) will remain in $\mathcal{K}(t)$
for all later times $t\in [t_0,T]$ provided it starts in
$\mathcal{K}(t_0)$ at $t_0$.

\hspace *{0.4cm}

In the following section we will give a proof of our theorem. In
forthcoming papers we will generalize the Hamilton-Ivey pinching
estimate and Hamilton's 4-dimensional theorem in [H2] to the case of
manifolds with boundary as applications of the various generalized
maximum principles.

\section {Proof of Theorem}

{\bf Proof}   \/  \/ We will adapt Hamilton's and  Chow and Lu's
idea (in [H2] and [CL]) to our case.

 We prove by contradiction. Suppose we have a solution $\sigma$ of the (PDE)
 which starts in $\mathcal{K}(t_0)$ at $t_0$ but runs out of $\mathcal{K}(t_2)$
  at some time $t_2\in (t_0,T]$. Then there exists
 $t_0 \le t_1 <t_2$ such that $\sigma(x,t_1)\in \mathcal{K}_x(t_1)$ for all $x\in M$ but for any $t\in (t_1,t_2]$
 there is some $x\in M$ such that $\sigma(x,t)\not\in \mathcal{K}_x(t)$. Let

 $g(x,t)=d(\sigma(x,t),\mathcal{K}_x(t))=inf\{|\sigma(x,t)-w|| w \in \mathcal{K}_x(t)\}$

 for each $x\in M, t\in [t_1,t_2]$, and

 $f(t)=sup_{x\in M}g(x,t)$

 for each $t\in [t_1,t_2]$, where we define the distant $d(w_1,w_2)$ between
 $w_1\in V_x$ and $w_2\in V_x$ using the metric $h_{ab}$ and denote it  by $|w_1-w_2|$ also.

By choosing $r$ large enough we may assume that $\sigma(x,t)\in
V(r)$ and  $d(\sigma(x,t),\partial (V(r)\cap
V_x))>d(\sigma(x,t),\mathcal{K}_x(t))$ for all $x\in M,t\in
[t_1,t_2]$, where $V(r)$ is the (closed) tubular neighborhood of the
zero section in $V$ whose intersection with each fiber $V_x$ is a
ball of radius $r$
 (measured by the bundle metric $h_{ab}$) around the origin. Note  that

$f(t)=sup\{l(\sigma(x,t)-v)|x\in M,v \in (\partial
\mathcal{K}_x(t))\cap V(r) \ \ and \ \ l\in S_v \mathcal{K}_x(t)\}$

for each $t\in (t_1,t_2)$.

We will prove by contradiction that $\frac{d^{+}f(t)}{dt}\le Cf(t)$
for any $t\in (t_1,t_2)$, where $\frac{d^{+}f(t)}{dt}=lim
sup_{h\rightarrow 0^+} \frac{f(t+h)-f(t)}{h}$, and $C=C_r$ is the
Lipschitz constant of $F(x,\sigma,t)$ (w.r.t. $\sigma$) within $
V(r)$. Suppose this is not the case.  Then we can find $t_a\in
(t_1,t_2)$ such that $\frac{d^{+}f}{dt}(t_a)> Cf(t_a)$.

From [CL] and the assumptions of our theorem we know that there
exist $x_\infty \in M, v_\infty \in (\partial
\mathcal{K}_{x_\infty}(t_a))\cap V(r)$, and $l_\infty \in
S_{v_\infty} \mathcal{K}_{x_\infty}(t_a)$, such that
$f(t_a)=l_\infty(\sigma(x_\infty,t_a)-v_\infty)=|\sigma(x_\infty,t_a)-v_\infty|$
and

$\frac{d^{+}f}{dt}(t_a)$

$\le
l_\infty(\Delta_{t_a}\sigma(x_\infty,t_a))+l_\infty(u^i(x_\infty,t_a)(\nabla_{t_a})_i
\sigma(x_\infty,t_a))+l_\infty(F(x_\infty,\sigma(x_\infty,t_a),t_a)-F(x_\infty,v_\infty,t_a))$

$\le
l_\infty(\Delta_{t_a}\sigma(x_\infty,t_a))+l_\infty(u^i(x_\infty,t_a)(\nabla_{t_a})_i
\sigma(x_\infty,t_a))+C|\sigma(x_\infty,t_a)-v_\infty|$.

Then we have

(*)\ \
$l_\infty(\Delta_{t_a}\sigma(x_\infty,t_a))+l_\infty(u^i(x_\infty,t_a)(\nabla_{t_a})_i
\sigma(x_\infty,t_a))>0$.

 Now we extend $v_\infty$ and
$l_\infty$ respectively by parallel translation ( defined by the
connection $\nabla_{t_a}$ ) along geodesics (w.r.t. $g_{ij}(t_a)$)
emanating radially from $x_\infty$, and we still denote what we get
by $v_\infty$ and $l_\infty$  respectively. Then by our assumption
on $\mathcal{K}$ we still have $v_\infty \in
\partial \mathcal{K}_x(t_a)$ and $l_\infty \in S_{v_\infty} \mathcal{K}_x(t_a)$. It is easy to see that the
function $l_\infty(\sigma(x,t_a)-v_\infty)$ has a local maximum at
$x_\infty$.

Then from (*) we have

$\Delta_{t_a}l_\infty(\sigma(x,t_a)-v_\infty)+u^i(x,t_a)(\nabla_{t_a})_i
l_\infty(\sigma(x,t_a)-v_\infty)>0$

at $x_\infty$ and hence also in a sufficiently small neighborhood of
$x_\infty$.
 It follows that $x_\infty\in \partial M$, and by Hopf's
lemma we have $\frac{\partial}{\partial
\nu}|_{x_\infty}l_\infty(\sigma(x,t_a)-v_\infty)
>0$. But on the other hand, by our assumption we have $v_\infty +(\nabla_{t_a})_\nu \sigma(x_\infty,t_a)\in
C_{v_\infty}\mathcal{K}_{x_\infty}(t_a)$ which implies
$l_\infty((\nabla_{t_a})_\nu \sigma(x_\infty,t_a))\le0$, so
$\frac{\partial}{\partial
\nu}|_{x_\infty}l_\infty(\sigma(x,t_a)-v_\infty)\le 0$, and we
arrive at a contradiction. Then by  Lemma 7 and Lemma 12 in [CL]
$f(t)=0$ on $(t_1,t_2]$, which contradicts the choice of $\sigma$,
and we are done.

\hspace *{0.4cm}

{\bf Remark } \ \ Similarly one can also generalize Chow and Lu's
Theorem 4 in [CL].

\hspace *{0.4cm}

\bibliographystyle{amsplain}

\hspace *{0.4cm}

{\bf Reference}

\bibliography{1}[BW] C. B$\ddot{o}$hm and B. Wilking, Manifolds with positive
 curvature operators are space forms, arXiv:math/0606187.

\bibliography{2}[BS] S. Brendle and R. M. Schoen, Manifolds with 1/4-pinched
 curvature are space forms, arXiv:math/0705.0766.

\bibliography{3}[CL] B. Chow and P. Lu, The maximum principle for systems of
parabolic equations subject to an avoidance set, Pacific J. Math.
214 (2004), 201-222.

\bibliography{4}[H1] R. S. Hamilton,Three-manifolds with positive Ricci curvature,
J. Diff. Geom. 17 (1982),255-306.

\bibliography{5}[H2]R. S. Hamilton, Four-manifolds with positive curvature curvature operator,
J. Diff. Geom. 24(1986),153-179.

\bibliography{6}[H3] R. S. Hamilton, The formation of singularities in the Ricci flow.
Surveys in differential geometry, Vol.II, 7-136, International
Press, 1995.

\bibliography{7}[I] T. Ivey, Ricci solitons on compact
three-manifolds, Diff. Geom. Appl. 3 (1993) 301-307.

\bibliography{8}[P] A. Pulemotov, Hopf boundary point lemma for vector bundle sections,
arXiv:math/0608037.

\bibliography{9}[S] Y. Shen, On Ricci deformation of a Riemannian metric on
manifolds with boundary, Pacific J. Math. 173 (1996) 203-221.

\end{document}